\documentclass[10pt,a4paper]{article}
\usepackage[a4paper, left=2cm, right=2cm, top=2cm]{geometry}

\usepackage[utf8]{inputenc}
\usepackage{amsmath, amsthm, amssymb, amsfonts}
\usepackage{siunitx}
\usepackage{natbib}
\usepackage[hidelinks]{hyperref}
\usepackage{doi}

\usepackage{subfigure}
\usepackage{authblk}
\usepackage{enumitem}
\usepackage{etoolbox}               

\usepackage{graphicx}
\usepackage{tikz}
\usepackage{pgfplots}

\usetikzlibrary{decorations.pathmorphing,arrows,intersections,arrows.meta,angles,quotes}

\pgfplotsset{
	kurze Legende/.style={
		legend image code/.code={
			\draw[##1,mark repeat=2,line width=0.6pt]
			plot coordinates {
				(0cm,0cm)
				(0.3cm,0cm)
			};
		}
	}
}

\pgfplotsset{
	compat = newest,
	scale only axis, 
	max space between ticks = 50pt,
	ticklabel style = {font=\footnotesize},
	legend style =  {font=\footnotesize},
	grid = major,
	grid style = {dotted},
	legend columns=1, 
	xtick pos=left,
	ytick pos=left
}

\pgfplotsset{select coords between index/.style 2 args={
		x filter/.code={
			\ifnum\coordindex<#1\fi
			\ifnum\coordindex>#2\fi
		}
}}

\definecolor{color1}{HTML}{0060AD} 
\definecolor{color2}{HTML}{FF4500} 
\definecolor{color3}{HTML}{FFA500} 
\definecolor{color4}{HTML}{006400} 
\definecolor{color5}{HTML}{9400D3} 
\definecolor{color6}{HTML}{800000} 
\definecolor{color7}{HTML}{000000} 
\definecolor{color8}{HTML}{0000FF} 
\definecolor{color9}{HTML}{FF0000} 
\definecolor{mycolor_blue}{RGB}{66,124,161}
\definecolor{mycolor_grey}{RGB}{198,198,198} 

\tikzstyle{line1} = [color=color7,semithick] 
\tikzstyle{line2} = [color=color2,semithick]
\tikzstyle{line3} = [color=color1,semithick]
\tikzstyle{line4} = [color=color5,semithick]
\tikzstyle{line5} = [color=color4,dash dot dot,semithick]
\tikzstyle{line6} = [color=color6,semithick]

\tikzstyle{mark1} = [color=color7,mark=x,mark size=2pt,mark options=solid,semithick] 
\tikzstyle{mark2} = [color=color2,mark=o,mark size=2pt,mark options=solid,semithick]
\tikzstyle{mark3} = [color=color1,mark=*,mark size=2pt,mark options=solid,semithick]
\tikzstyle{mark4} = [color=color5,mark=triangle,mark size=2pt,mark options=solid,semithick]
\tikzstyle{mark5} = [color=color4,mark=square,mark size=2pt,mark options=solid,semithick]
\tikzstyle{mark6} = [color=color7,mark=o,mark size=2pt,mark options=solid,semithick]
\tikzstyle{mark7} = [color=color7,mark=*,mark size=2pt,mark options=solid,semithick]
\tikzstyle{mark8} = [color=color7,mark=triangle,mark size=2pt,mark options=solid,semithick]

\usetikzlibrary{external}
\title{Discrete Adjoint Momentum-Weighted Interpolation Strategies}
\author[]{Niklas K\"uhl\thanks{niklas.kuehl@tuhh.de} }
\author[]{Thomas Rung}

\affil[]{Hamburg University of Technology, Institute for Fluid Dynamics and Ship Theory, Am Schwarzenberg-Campus 4, D-21075 Hamburg, Germany}

\begin{document}

\providetoggle{tikzExternal}
\settoggle{tikzExternal}{false}

\maketitle

%
\section{Introduction}
This Technical Note outlines an adjoint complement to a critical building block of pressure-based Finite-Volume (FV) flow solvers that employ a collocated variable arrangement to simulate virtually incompressible fluids, cf. \cite{ferziger2012computational}. The focal point is to strengthen the adjoint pressure-velocity coupling by using an adjoint Momentum-Weighted Interpolation (MWI) strategy. To this end, analogies of established primal MWI techniques are derived and investigated. The study reveals the merits of an adjoint MWI but also highlights the importance of its careful implementation.

The interplay between the pressure and momentum sources is frequently addressed by MWI methods, which support the prediction of mass fluxes through control volume faces in combination with pressure-based FV methods and collocated variable arrangements, cf. \cite{mencinger2007finite, choi2003use}. MWI aims at stabilizing the algorithm by suppressing an odd-even decoupling between the pressure and the velocities. Based on the seminal work of \cite{rhie1983numerical}, such methods typically gather under the umbrella of Rhie-Chow Interpolation (RCI) schemes and represent a standard in FV solvers dedicated to industrial problems, cf. \cite{yu2002checkerboard, pascau2011cell, bartholomew2018unified}. The technique is of discrete nature, disappears in a continuous framework, and can thus form a bridge between continuous derivation and discrete adjoint implementation, cf. \cite{stuck2013adjoint, kroger2018adjoint, kuhl2021adjoint}.

When attention is drawn to an adjoint optimization framework, the interest in an adjoint MWI follows from the increased number of source terms of the adjoint momentum equation compared to its primal companion. Any variation of the primal fluid velocity -- which is not linked to the adjoint velocity -- results in an adjoint source term, which is mostly treated explicitly. Examples include velocity depending volume cost functional or the convective transport of auxiliary variables, e.g., turbulence parameters. Consistent and numerically robust treatment of the many explicit source terms is desirable and will be assessed in this note. In the following, vectors and tensors are defined with reference to Cartesian coordinates and denoted by Latin subscripts. Einstein’s summation is used for Latin subscripts.

%
\section{Primal Flow and Starting Point}
The governing flow equations for the primal velocity $v_\mathrm{i}$ and pressure $p$ of an incompressible viscous ($\mu$) fluid with density $\rho$ yield the following residual $r^\mathrm{(\cdot)}$ form
\begin{alignat}{2}
r_\mathrm{i}^\mathrm{v_\mathrm{i}} &= 0 &&= v_\mathrm{k} \rho \frac{\partial v_\mathrm{i}}{\partial x_\mathrm{k}}  - \frac{\partial }{\partial x_\mathrm{k}} \left[ \mu \left( \frac{\partial v_\mathrm{i}}{\partial x_\mathrm{k}} + \frac{\partial v_\mathrm{k}}{\partial x_\mathrm{i}} \right) \right] + \frac{\partial p}{\partial x_\mathrm{i}} - q_\mathrm{i} \, , \label{equ:primal_momentum_balance} \\
r^\mathrm{p} &= 0 &&= -\frac{\partial v_\mathrm{k}}{\partial x_\mathrm{k}} \label{equ:primal_mass_balance} \, .
\end{alignat}
Here $q_\mathrm{i}$ refers to a general differentiable momentum source, e.g., $q_\mathrm{i} = \rho g_\mathrm{i}$ in case of an acceleration by $g_\mathrm{i}$. Since the note is devoted to the numerical treatment of adjoint momentum sources, which often follow from transport equations of auxiliary variables, we define a generic stationary transport equation for the auxiliary quantity $\varphi$ by the following residual form
\begin{align}
r^\mathrm{\varphi} = 0 = v_\mathrm{k} \frac{\partial \varphi}{\partial x_\mathrm{k}} - \frac{\partial}{\partial x_\mathrm{k}} \left[ \mu^\varphi \frac{\partial \varphi}{\partial x_\mathrm{k}} \right] - s^\mathrm{\varphi} \label{equ:primal_generic_balance} \, .
\end{align}
Here, $\mu^\varphi$ and $s^\mathrm{\varphi}$ refer to a $\varphi$-specific diffusivity and source term, respectively. Typical auxiliary equations address the statistical modeling of turbulence, e.g. $\varphi \to k$ and $\varphi \to \varepsilon$ in case of the classical $k-\varepsilon$ turbulence model, cf. \cite{jones1972prediction}, or the transport of a volume fraction $\varphi \to c$ to model two-phase flow, e.g. along the route of the Volume of Fluid (VoF) method proposed by \cite{hirt1981volume}.
\subsection{Starting Point of the Momentum-Weighted Interpolation}
\begin{figure}[!h]
\centering
\iftoggle{tikzExternal}{
\input{./tikz/finite_volume_approximation_field.tikz}
}{
\includegraphics{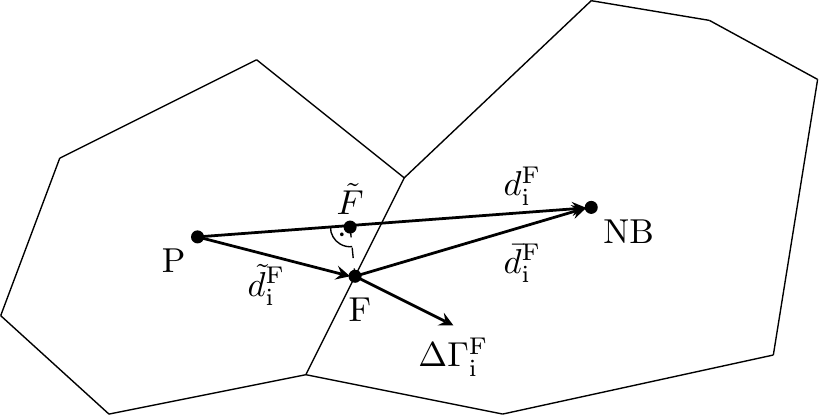}
}
\caption{Schematic representation of a cell-centered Finite-Volume arrangement for the cell around a focal center P, a neighbor cell around the center location NB, and their shared face $\Delta \Gamma_\mathrm{i}^\mathrm{F}$.}
\label{fig:finite_volume_approximation}
\end{figure}

In a situation at rest ($v_\mathrm{i} \to 0$) the momentum balance equation (\ref{equ:primal_momentum_balance}) simplifies to a ''hydrostatic'' relation, viz. 
\begin{align}
\mathrm{C_\mathrm{1}:} \quad \frac{\partial p}{\partial x_\mathrm{i}} = q_\mathrm{i}
\qquad \qquad \mathrm{or} \qquad \qquad
\mathrm{C_\mathrm{2}:} \quad \frac{\partial p}{\partial x_\mathrm{i}} =
\frac{\partial \, q_\mathrm{k} x_\mathrm{k}}{\partial x_\mathrm{i}} - \frac{\partial q_\mathrm{k}}{\partial x_\mathrm{i}} x_\mathrm{k}
\, . \label{equ:primal_momentum_source}
\end{align}
Equation (\ref{equ:primal_momentum_source}) distinguishes two -- continuously speaking identical -- formulation concepts ($\mathrm{C_\mathrm{1}},\mathrm{C_\mathrm{2}}$) for the momentum source $q_\mathrm{i}$, which provide access to a -- discretely speaking -- robust approximation, cf.  \cite{bartholomew2018unified}. Using a cell-centered storage approach and a conventional second-order accurate mid-point integration rule, a FV approximation of the spatial integral over the two alternatives reads
\begin{align}
    \mathrm{C_\mathrm{1}:} \, \int_{\Delta \Omega (P)} q_\mathrm{i}
    \mathrm{d} \Omega 
    \approx \left[ q_\mathrm{i} \Delta \Omega \right]^\mathrm{P} \, ,
    \quad
    \mathrm{C_\mathrm{2}:} \, \int_{\Delta \Omega (P)} q_\mathrm{i} \mathrm{d} \Omega
    \approx
    \oint_{\Gamma(P)} q_\mathrm{k} x_\mathrm{k} \mathrm{d} \Gamma_\mathrm{i} - x_\mathrm{k}^\mathrm{P} \oint_{\Gamma(P)} q_\mathrm{k} \mathrm{d} \Gamma_\mathrm{i}
    \approx    
    \sum_{\Delta \Gamma(P)} \big[ x_\mathrm{k}^\mathrm{F} - x_\mathrm{k}^\mathrm{P} \big] q_\mathrm{k}^\mathrm{F} \Delta \Gamma_\mathrm{i}^\mathrm{F} \; . 
\end{align}
Here $P$ denotes the center of the control volume (CV), $\Delta \Omega$ refers to its size, and $F$ marks the face center locations of the discrete, outward-pointing face segment vectors $\Delta \Gamma_\mathrm{i}$ enclosing the CV, cf. Fig \ref{fig:finite_volume_approximation}. Later on, MWI strategies employ volume averaged quantities that read
\begin{align}
    Q_\mathrm{i}^{\mathrm C_\mathrm{1}} = \frac{1}{\Delta \Omega(P)} \int_{\Delta \Omega (P)} q_\mathrm{i} \mathrm{d} \Omega 
    \approx q_\mathrm{i}^\mathrm{P}
    \qquad \qquad \mathrm{and} \qquad \qquad
    Q_\mathrm{i}^{\mathrm C_\mathrm{2}} = \frac{1}{\Delta \Omega (P)}\sum_{\Delta \Gamma(P)} \big[ x_\mathrm{k}^\mathrm{F} - x_\mathrm{k}^\mathrm{P} \big] q_\mathrm{k}^\mathrm{F} \Delta \Gamma_\mathrm{i}^\mathrm{F}  \, .
\label{eq:meanQ}
\end{align}
A crucial aspect refers to the interpolation of $q_\mathrm{i}^\mathrm{F}$ to the face centers, as seen in Sec. \ref{sec:application}.

\section{Adjoint Flow}
Adjoint procedures are an efficient approach to identify objective functional sensitivities with respect to (w.r.t.) a large number of optimization parameters, cf. \cite{jameson1995optimum}. A general integral objective functional $J$ which consists of boundary ($j^\mathrm{\Gamma}$) and volume ($j^\mathrm{\Omega}$) contributions along the objective surface $\Gamma^\mathrm{O}$ and the objective volume $\Omega^\mathrm{O}$ reads
\begin{align}
J = \int_{\Gamma^\mathrm{O}} j^\mathrm{\Gamma} \mathrm{d} \Gamma + \int_{\Omega^\mathrm{O}} j^\mathrm{\Omega} \mathrm{d} \Omega \label{equ:general_objective} \, . 
\end{align}
This is supplemented by the primal residuals (\ref{equ:primal_momentum_balance}), (\ref{equ:primal_mass_balance}) and (\ref{equ:primal_generic_balance}) to form an augmented objective functional, often referred to as a Lagrange functional
\begin{align}
L = J  + \int_\Omega ( \hat{v}_\mathrm{i} \, r_\mathrm{i}^\mathrm{v_\mathrm{i}} + \hat{p} \, r^\mathrm{p} + \hat{\varphi} \, r^\mathrm{\varphi}) \, \mathrm{d} \Omega \, .
\end{align}
Quantities in front of the residual expressions denote the adjoint velocity $\hat{v}_\mathrm{i}$, the adjoint pressure $\hat{p}$, and the adjoint auxiliary quantity $\hat{\varphi}$. The differentiation the Lagrangian in the direction of the primal state is introduced based on first-order optimality conditions, i.e. $\delta_{v_\mathrm{i}} L \cdot \delta v_\mathrm{i} \overset{!}{=} 0 \, \forall \, \delta v_\mathrm{i}$, $\delta_{p} L \cdot \delta p \overset{! }{=} 0 \, \forall \, \delta p$, $\delta_{\varphi} L \cdot \delta \varphi \overset{!}{=} 0 \, \forall \, \delta \varphi$, to derive the governing adjoint equations, cf. \cite{lohner2003adjoint, othmer2008continuous, giannakoglou2008adjoint, stuck2013adjoint, kroger2018adjoint, kuhl2019decoupling, kuhl2021adjoint, kuhl2021adjoint_2, kuhl2021continuous}, viz. 
\begin{alignat}{2}
r_\mathrm{i}^\mathrm{\hat{v}_\mathrm{i}} &= 0 &&= - v_\mathrm{k} \rho \frac{\partial \hat{v}_\mathrm{i}}{\partial x_\mathrm{k}} - \frac{\partial }{\partial x_\mathrm{k}} \left[ \mu \left( \frac{\partial \hat{v}_\mathrm{i}}{\partial x_\mathrm{k}} + \frac{\partial \hat{v}_\mathrm{k}}{\partial x_\mathrm{i}} \right) \right] + \frac{\partial \hat{p}}{\partial x_\mathrm{i}} - \hat{v}_\mathrm{i} \frac{\partial q_\mathrm{k}}{\partial v_\mathrm{k}} - \hat{\varphi} \frac{\partial s^\mathrm{\varphi}}{\partial v_\mathrm{i}} + \hat{v}_\mathrm{k} \rho \frac{\partial v_\mathrm{k}}{\partial x_\mathrm{i}} + \hat{\varphi} \frac{\partial \varphi}{\partial x_\mathrm{i}} + \frac{\partial j^\mathrm{\Omega}}{\partial v_\mathrm{i}} \label{equ:adjoint_momentum_balance} \, , \\
r^\mathrm{\hat{p}} &= 0 &&= -\frac{\partial \hat{v}_\mathrm{k}}{\partial x_\mathrm{k}} - \hat{v}_\mathrm{i} \frac{\partial q_\mathrm{i}}{\partial p} - \hat{\varphi} \frac{\partial s^\mathrm{\varphi}}{\partial p} + \frac{\partial j^\mathrm{\Omega}}{\partial p} \label{equ:adjoint_mass_balance} \, , \\
r^\mathrm{\hat{\varphi}} &= 0 &&= -v_\mathrm{k} \frac{\partial \hat{\varphi}}{\partial x_\mathrm{k}} - \frac{\partial}{\partial x_\mathrm{k}} \left[ \mu^\mathrm{\hat{\varphi}} \frac{\partial \hat{\varphi}}{\partial x_\mathrm{k}} \right] - \hat{v}_\mathrm{i} \frac{\partial q_\mathrm{i}}{\partial \varphi} - \hat{\varphi} \frac{\partial s^\mathrm{\varphi}}{\partial \varphi} + \frac{\partial j^\mathrm{\Omega}}{\partial \varphi} \label{equ:adjoint_generic_balance} \, .
\end{alignat}
The adjoint equations (\ref{equ:adjoint_momentum_balance})-(\ref{equ:adjoint_generic_balance}) display significantly more source terms than their primal companions (\ref{equ:primal_momentum_balance})-(\ref{equ:primal_generic_balance}). Depending on the source term magnitude, this might be obstructive for the adjoint pressure-velocity coupling. Remaining optimality criteria and a correct approximation of the primal and adjoint equations allow for a sensitivity (\cite{othmer2008continuous, stuck2013adjoint, kuhl2021adjoint_2}) rule along the controlled design wall
\begin{align}
    \delta_\mathrm{u} J = \int_\mathrm{\Gamma^\mathrm{D}} s \mathrm{d} \Gamma
    \qquad \qquad \mathrm{with} \qquad \qquad
    s = - \mu \frac{\partial v_\mathrm{i}}{\partial n} \frac{\partial \hat{v}_\mathrm{i}}{\partial n} \, . \label{equ:sensitivity_derivative}
\end{align}
Note that the locations of the control ($\Gamma^\mathrm{D}$) and the objective ($\Gamma^\mathrm{O}$) do not necessarily coincide, cf. \cite{kuhl2019decoupling}. Analogous to the primal relation (\ref{equ:primal_momentum_source}), the adjoint momentum balance simplifies for a vanishing adjoint velocity ($\hat{v}_\mathrm{i} \to 0$) towards an "adjoint hydrostatic" pressure balance, viz.
\begin{align}
\mathrm{C_\mathrm{1}:} \quad \frac{\partial \hat{p}}{\partial x_\mathrm{i}} = 
 \hat q_\mathrm{i}
\qquad \mathrm{or} \qquad
\mathrm{C_\mathrm{2}:} \quad \frac{\hat \partial p}{\partial x_\mathrm{i}} =
\frac{\partial \, \hat q_\mathrm{k} x_\mathrm{k}}{\partial x_\mathrm{i}} - \frac{\partial \hat q_\mathrm{k}}{\partial x_\mathrm{i}} x_\mathrm{k} \, ,
\qquad \mathrm{with}  \qquad
\hat{q}_\mathrm{k} = \hat{\varphi} \frac{\partial s^\mathrm{\varphi}}{\partial v_\mathrm{i}} +  \hat{\varphi} \frac{\partial \varphi}{\partial x_\mathrm{i}} + \frac{\partial j^\mathrm{\Omega}}{\partial v_\mathrm{i}} \, . \label{equ:adjoint_momentum_source}
\end{align}
Again two alternative formulation concepts ($C_\mathrm{1}, C_\mathrm{2}$) exist, which might help to improve the adjoint coupling. Compared to the primal system, the adjoint momentum sources tend to be much more volatile with potentially more significant local gradients, e.g., in free surface models that feature a  ''jump'' of the mixture fraction $c$ across the compressive interface ($\nabla_\mathrm{i} \varphi \to \nabla_\mathrm{i} c$). Hence, a robust and, at the same time, consistent adjoint pressure-velocity coupling is of significance and will be addressed in the following section.

%
\section{Discrete Adjoint Momentum-Weighted Interpolation}
The discrete implementation of primal pressure-velocity coupling strategies has been the subject of various research efforts and will not be revisited here. A recent overview of modern MWI methods and their implementation is given in \cite{bartholomew2018unified}. Aiming to adopt these concepts in an adjoint context, we start from the semi-discrete version of the adjoint momentum balance (\ref{equ:adjoint_momentum_balance}) for a CV around P based on a cell-centered FV discretization, cf. Fig. \ref{fig:finite_volume_approximation}
\begin{align}
A^\mathrm{P} \hat{v}_\mathrm{i}^\mathrm{P} + \sum_{NB(P)} a^\mathrm{NB} \hat{v}_\mathrm{i}^\mathrm{NB} = \left[ \hat{q}_\mathrm{i} - \frac{\partial \hat{p}}{\partial x_\mathrm{i}} \right]^\mathrm{P} \Delta \Omega^\mathrm{P} \label{equ:semi_discrete_adjoint_momentum} \, .
\end{align}
Here $A^\mathrm{P}$ and $a^\mathrm{NB}$ denote the main and off diagonal coefficients, respectively. A similar expression is obtained for each neighboring cell center velocity $\hat v_\mathrm{i}^{\mathrm{NB}}$. Similar to the discretized primal continuity equation (\ref{equ:primal_mass_balance}), the discretized adjoint continuity equation (\ref{equ:adjoint_mass_balance}) requires face-based velocities $\hat v_\mathrm{i}^{\mathrm F}$ to compute fluxes, viz.
\begin{align}
    \sum_{\Delta \Gamma(P)} \hat v_\mathrm{i}^{\mathrm F} \Delta \Gamma_\mathrm{i}^\mathrm{F} = -\left[  \hat{v}_\mathrm{k} \frac{\partial q_\mathrm{k}}{\partial p} + \hat{\varphi} \frac{\partial s^\mathrm{\varphi}}{\partial p} - \frac{\partial j^\mathrm{\Omega}}{\partial p}\right]^{\mathrm P} \Delta \Omega^{\mathrm P}
\end{align}
The derivation follows analogously to primal RCI strategies (\cite{rhie1983numerical, pascau2011cell, bartholomew2018unified}) and yields the following adjoint MWI rule
\begin{align}
\hat{v}_\mathrm{i}^\mathrm{F}  = \overline{\hat{v}_\mathrm{i}^\mathrm{F}} + \overline{\frac{\partial \hat{v}_\mathrm{i}}{\partial x_\mathrm{k}} \bigg|^\mathrm{F}} \left( x_\mathrm{k}^\mathrm{F} - x_\mathrm{k}^\mathrm{\tilde{F}} 
\right)- \beta^\mathrm{F} \left[ \left( \frac{\partial \hat{p}}{\partial x_\mathrm{i}}  - \hat{Q}_\mathrm{i} \right)^\mathrm{C1} - \left( \frac{\partial \hat{p}}{\partial x_\mathrm{i}} - \hat{Q}_\mathrm{i}\right)^\mathrm{C2} \right]^\mathrm{F} \, , \label{equ:adjoint_rhie_chow}
\end{align}
where over-lined expressions $\overline{(\cdot)^\mathrm{F}} = \lambda^\mathrm{F} (\cdot)^\mathrm{NB} + (1 - \lambda^\mathrm{F}) (\cdot)^\mathrm{P}$ follow from a simple linear interpolation into the perpendicular location $\tilde F$ using the geometric weight $\lambda^\mathrm{F} = (x_\mathrm{k}^\mathrm{F} - x_\mathrm{k}^\mathrm{P}) d_\mathrm{k}^\mathrm{F} / (d_\mathrm{m}^\mathrm{F} d_\mathrm{m}^\mathrm{F})$. The second right-hand side (r.h.s.) term refers to an optional non-orthogonality correction, cf. Fig. \ref{fig:finite_volume_approximation}. 
The four remaining contributions on the r.h.s. comprise the adjoint MWI rule. The coefficient reads $\beta^\mathrm F$ reads $(d_\mathrm{k} \Delta \Gamma_\mathrm{k}^\mathrm F)/ \, \overline{ A^\mathrm F}$. The pressure contributions are computed from 
\begin{align}
    \frac{\partial \hat{p}}{\partial x_\mathrm{i}}\bigg|^\mathrm{F,C1} = \frac{\hat{p}^\mathrm{NB} - \hat{p}^\mathrm{P}}{||d_\mathrm{k}^\mathrm{F}||} \frac{d_\mathrm{i}^\mathrm{F}}{||d_\mathrm{k}^\mathrm{F}||} 
    \qquad \qquad \text{and} \qquad \qquad
    \frac{\partial \hat{p}}{\partial x_\mathrm{i}}\bigg|^\mathrm{F,C2} = \overline{\frac{\partial \hat{p}}{\partial x_\mathrm{i}}\bigg|^\mathrm{F}} 
    = \lambda^\mathrm{F} \frac{\partial \hat{p}}{\partial x_\mathrm{i}}\bigg|^\mathrm{P} +  \left( 1 - \lambda^\mathrm{F} \right)  \frac{\partial \hat{p}}{\partial x_\mathrm{i}}\bigg|^\mathrm{NB}
    \label{equ:discrete_gradadp_C1_C2} \, ,
\end{align}
and essentially correspond to their classical RCI companion (\cite{rhie1983numerical}), cf. \cite{bartholomew2018unified}. The source term expressions $\hat Q_\mathrm{i}^{\mathrm{C1}}$ and $\hat Q_\mathrm{i}^{\mathrm{C2}}$  follow from (\ref{eq:meanQ}, \ref{equ:adjoint_momentum_source}) and are reconstructed at the face F using
\begin{align}
    \hat{Q}_\mathrm{i}^\mathrm{C1, F} = \overline{\overline{\hat{Q}_\mathrm{i}^\mathrm{C1, F}}} = \left( 1 - \lambda^\mathrm{F} \right) \hat{Q}_\mathrm{i}^\mathrm{C1, P} + \lambda^\mathrm{F} \hat{Q}_\mathrm{i}^\mathrm{C1, NB}
    \qquad \text{and} \qquad
    \hat{Q}_\mathrm{i}^\mathrm{C2, F} = \overline{\hat{Q}_\mathrm{i}^\mathrm{C2, F}} = \lambda^\mathrm{F} \hat{Q}_\mathrm{i}^\mathrm{P, C2} + \left( 1 - \lambda^\mathrm{F} \right) \hat{Q}_\mathrm{i}^\mathrm{NB, C2} \label{equ:discrete_adq_C1_C2} \, .
\end{align}
The cell-centered evaluation of $\hat Q_\mathrm{i}$ adheres to different approximations, cf. (\ref{eq:meanQ}), viz.
\begin{align}
    \hat{Q}_\mathrm{i}^\mathrm{C1, P} \approx  \hat{q}_\mathrm{i}^\mathrm{P} \qquad {\mathrm{and}} \qquad
 \hat{Q}_\mathrm{i}^\mathrm{C2, P}   \approx \frac{1}{\Delta \Omega^\mathrm{P}} \sum_{F(P)} \big[ x_\mathrm{k}^\mathrm{F} - x_\mathrm{k}^\mathrm{P} \big] \overline{\overline{\hat{q}_\mathrm{k}^\mathrm{F}}} \Delta \Gamma_\mathrm{i}^\mathrm{F} \, 
. \label{equ:discrete_adbf_C2}
\end{align}
A crucial aspect refers to the face reconstructions of adjoint momentum source contributions with  reversed levers $\overline{\overline{(\cdot)^\mathrm{F}}} = (1 - \lambda^\mathrm{F}) (\cdot)^\mathrm{NB} + \lambda^\mathrm{F} (\cdot)^\mathrm{P}$. This follows from primal MWI strategies and stems from the fact that the influence of volume forces increases with greater leverage (\cite{bartholomew2018unified, mencinger2007finite}). 

%
\section{Application}
\label{sec:application}
This section assesses the credibility of different adjoint MWI approaches. The investigations refer to the laminar two-phase flow around a two-dimensional submerged circular cylinder at rest and involves a surface-based drag cost functional, i.e. $j^\mathrm{\Gamma} = (p \delta_\mathrm{ik} - \mu (\nabla_\mathrm{k} v_\mathrm{i} + \nabla_\mathrm{i} v_\mathrm{k})) n_\mathrm{k} \delta_\mathrm{i1}$ -- where $n_\mathrm{k}$ and $\delta_\mathrm{ik}$ denote the local boundary normal and the Kronecker-Delta.

\begin{figure}[!ht]
\centering
\subfigure[]{
\iftoggle{tikzExternal}{
\input{./tikz/2D/cylinder_scetch.tikz}
}{
\includegraphics{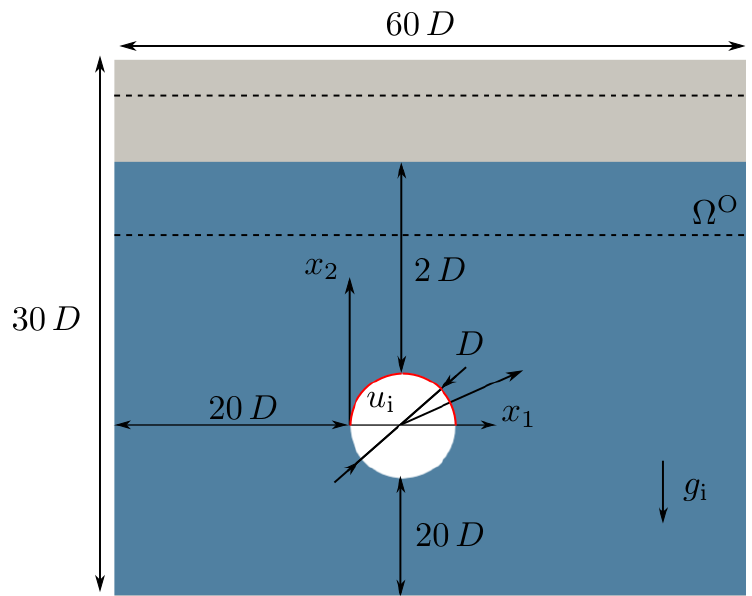}
}
}
\subfigure[]{
\includegraphics[scale=1]{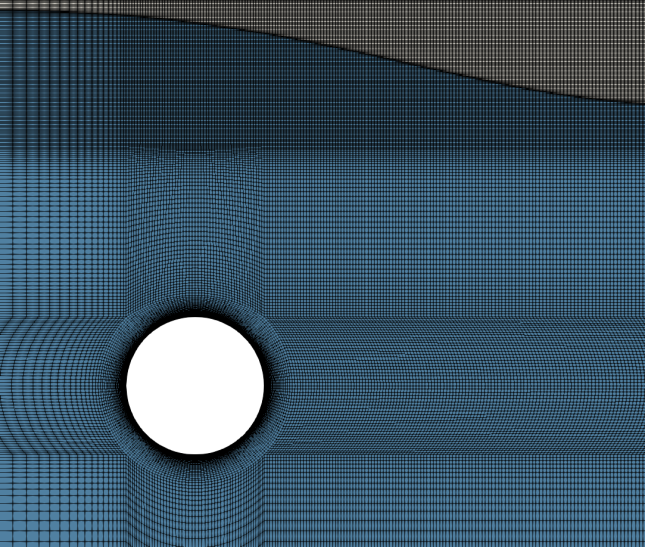}
}
\caption{Submerged cylinder case ($\mathrm{Re}_\mathrm{D} = 20$, $\mathrm{Fn}=0.75$); (a) Schematic  of the initial configuration around the controlled cylinder shape $u_\mathrm{i}$ (red) and (b) portion of the structured numerical grid around the cylinder.}
\label{fig:cylinder_fn_075}
\end{figure}

As illustrated in Fig. \ref{fig:cylinder_fn_075} (a), the origin of the cylinder is positioned two and a half diameters $D$ underneath an initial calm-water free surface. The employed two-dimensional domain features a length and height of $60 \, D$ and $30 \, D$, where the inlet and bottom boundaries are located 20,5 diameters away from the cylinder's origin. At the inlet, a homogeneous unidirectional  (horizontal) bulk flow $v_\mathrm{i} = v_\mathrm{1} \delta_{i1}$ is imposed for both phases in conjunction with a calm water concentration distribution. Slip walls are used along the top and bottom boundaries, and a hydrostatic pressure boundary is employed along with the outlet. The grid is stretched in the longitudinal direction ($x_\mathrm{1}$) towards the outlet to suppress the outlet wave field and comply with the outlet condition. The study is performed at $\mathrm{Re}_\mathrm{D} = v_\mathrm{1} D/\nu^\mathrm{ b} = \SI{20}{}$ and $\mathrm{Fn} = v_\mathrm{1}/\sqrt{G \, 2 \, D} = \SI{0.75}{}$, based on the gravitational acceleration $G$, the inflow velocity $v_\mathrm{1}$ and the kinematic viscosity of the water $\nu^\mathrm{ b}$. The expected dimensionless wave length reads $\lambda/ D = 2 \, \pi \, \mathrm{Fn}^2 = 3.534$. To ensure the independence of the objective functional value w.r.t. spatial discretization, a grid study was conducted prior to the optimization study. Part of the utilized structured numerical grid is displayed in Fig. \ref{fig:cylinder_fn_075} (b). It consists of approximately $\SI{215000}{}$ control volumes where the cylinder shape is discretized with 500 surface elements along the circumference. The non-dimensional wall-normal distance of the first grid layer reads $y^+ \approx \SI{0.01}{}$ and the refined grid in the free surface region employs isotropic spacing with $\Delta x_\mathrm{1} =\Delta x_\mathrm{2} \approx \lambda/100$.

A pressure-based, second-order accurate FV scheme using a cell-centered, co-located variable arrangement is employed to approximate the partial differential equations of the primal and adjoint systems, cf. \cite{rung2009challenges, stuck2012adjoint, kroger2016numerical, kuhl2021phd}. The sequential procedure employs a SIMPLE-type pressure correction method, cf. \cite{yakubov2015experience}, and the solution is iterated to convergence until all -- or only a selected subset of equations -- are converged below a prescribed global tolerance $\mathrm{R}^\mathrm{\varphi,min}$. Convective primal [adjoint] momentum fluxes are approximated using the QUICK [QDICK] scheme, cf. \cite{stuck2013adjoint}. The approximation of the concentration equation has been outlined in \cite{kuhl2021adjoint}, where the traditional VoF approach follows from a  compressive primal/hybridized continuous-discrete adjoint HRIC scheme. Using an Euler implicit approach, the simulations are advanced to a steady state in pseudo time.

\begin{figure}[h]
\centering
\iftoggle{tikzExternal}{
\input{./tikz/2D/cylinder_fd_drift_inverse.tikz}
}{
\includegraphics{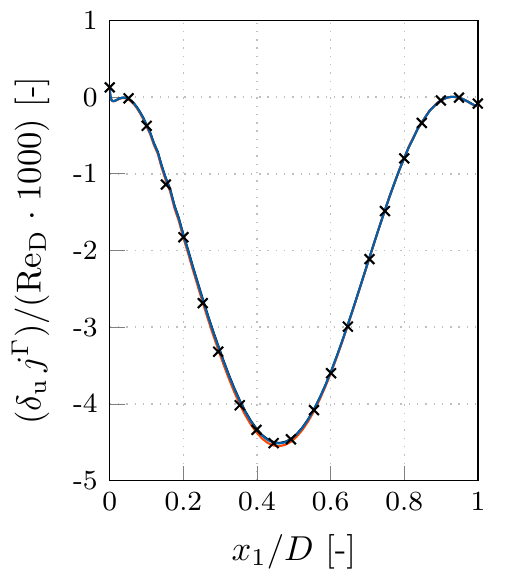}
\includegraphics{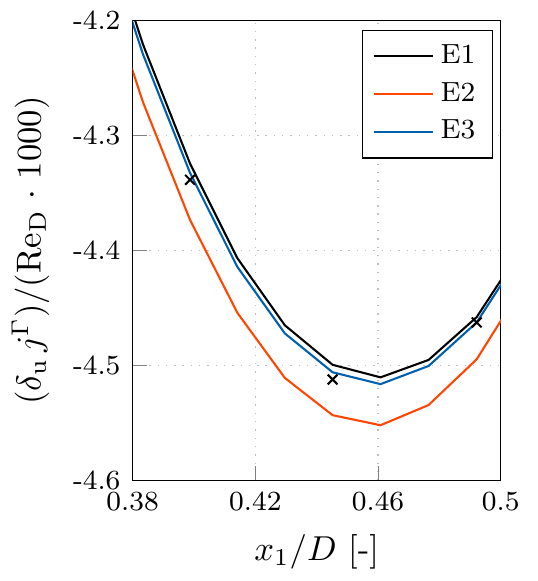}
\includegraphics{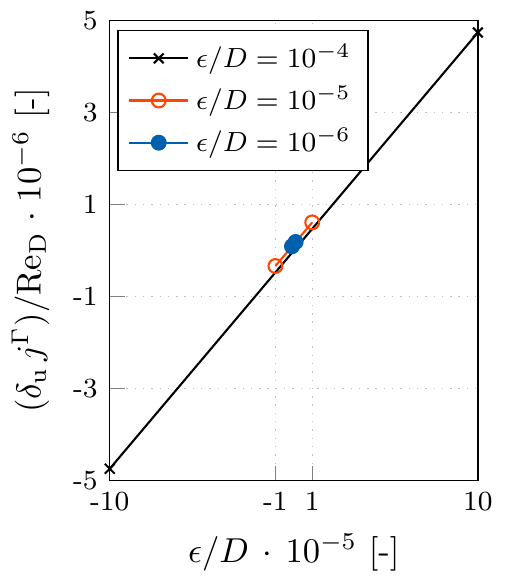}
}
\caption{Submerged cylinder case ($\mathrm{Re}_\mathrm{D} = 20$, $\mathrm{Fn}=0.75$);  Adjoint and Finite-Difference (FD) based sensitivity derivative for a drag functional along  the  upper cylinder (left) and magnification of the area of greatest sensitivity (center), supplemented by three exemplary finite system answers at $x_\mathrm{1} / D = 1/4$ (right).}
\label{fig:cylinder_fd_study}
\end{figure}

Emphasis is put on the predictive agreement of the adjoint sensitivities with results of a Finite Difference (FD) approach and the iterative convergence behavior. Mind that the adjoint momentum balance experiences an additional source that follows from the primal convective concentration transport which in turn features huge local gradients, i.e. $\hat{q}_\mathrm{i} \gets \hat{c} \nabla_\mathrm{i} c$, cf.  (\ref{equ:adjoint_momentum_source}). Three numerical experiments E1-E3 are compared, which differently include the adjoint momentum source into the adjoint MWI procedure (\ref{equ:adjoint_rhie_chow}):
\begin{itemize}
    \item[E1:] The first experiment neglects all adjoint momentum sources $\hat Q_\mathrm{i}$ within the adjoint MWI (\ref{equ:adjoint_rhie_chow}) and considers only the adjoint pressure gradients in line with (\ref{equ:discrete_gradadp_C1_C2}).
    \item[E2:] A second experiment includes the adjoint volume forces  $\hat Q_\mathrm{i}$  within the adjoint MWI (\ref{equ:adjoint_rhie_chow}), but invariably uses a classical linear interpolation, i.e. $\overline{\overline{(\cdot)^\mathrm{F}}} \to \overline{(\cdot)^\mathrm{F}}$ in Eqns. (\ref{equ:discrete_adq_C1_C2})-(\ref{equ:discrete_adbf_C2}).
    \item[E3:] Only the third approach employs twisted levers for $\overline{\overline{(\cdot)^\mathrm{F}}} $ and computes the adjoint face velocities $\hat{v}_\mathrm{i}^\mathrm{F}$ as described in (\ref{equ:adjoint_rhie_chow})-(\ref{equ:discrete_adbf_C2}).
\end{itemize}
The comparison between adjoint sensitivities and Finite-Differences (FD) is depicted in Fig. \ref{fig:cylinder_fd_study}. The credibility of the FD results is ensured by verifying the linearity of the FD-analysis using three perturbation magnitudes  $\epsilon / D \in [10^{-4}, 10^{-5}, 10^{-6} ] $. The  control is  restricted to the upper half of the cylinder (cf. Fig \ref{fig:cylinder_fn_075}a), for which FD results are extracted at 21 discrete positions. To this end, 42 additional simulations were performed to obtain second-order accurate central differences. An exemplary documentation of the systems linear answer is displayed in the right graph of Fig. \ref{fig:cylinder_fd_study}, which refers to an exemplary surface location $x_\mathrm{1} / D = 1/4$.

\begin{figure}[h]
\centering
\iftoggle{tikzExternal}{
\input{./tikz/2D/cylinder_adjoint_residuals_drag.tikz}
}{
\includegraphics{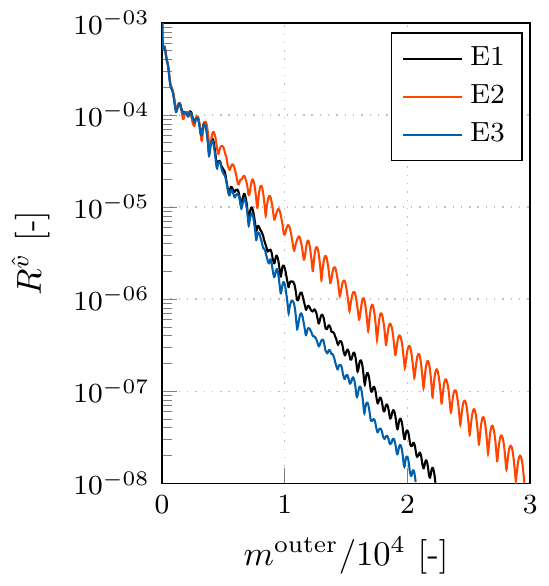}
\includegraphics{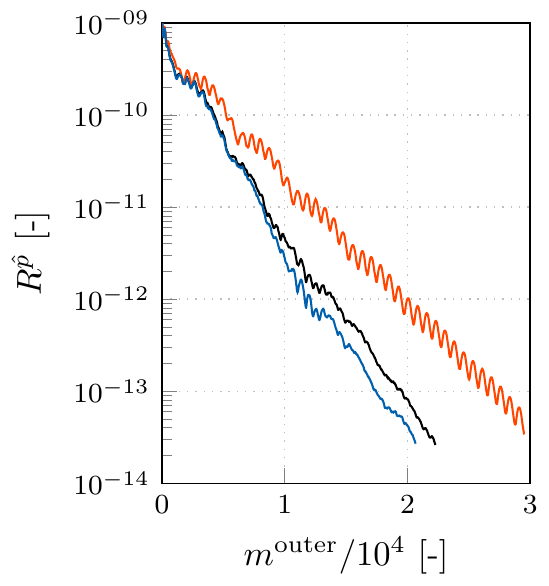}
\includegraphics{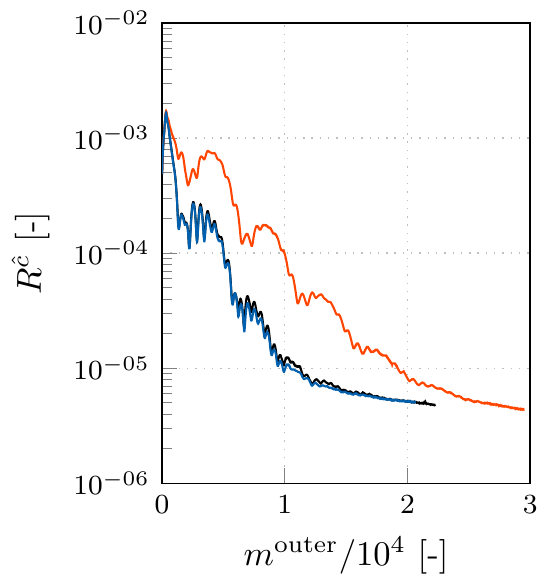}
}
\caption{Submerged cylinder case ($\mathrm{Re}_\mathrm{D} = 20$, $\mathrm{Fn}=0.75$); Global adjoint residuals (\ref{equ:global_residual}) over the outer iterations of (left) the adjoint velocity -- summed up for both spatial directions --, (center) adjoint pressure, and (right) adjoint concentration equation for the investigated adjoint MWI configurations E1-E3. All simulations terminated after an adjoint velocity residual value of $R^{\hat{v}} \le 10^{-08}$.}
\label{fig:cylinder_adjoint_residuals_drag}
\end{figure}

All adjoint sensitivities E1-E3 essentially agree fairly with the FD results. However, a closer inspection of the maximum sensitivities (center) reveals that the adjoint sensitivities do not coincide and feature minor deviations. The most significant deviations are observed in conjunction with the wrong levers (E2). The prediction improves if adjoint momentum sources are completely neglected (E1), but a superior agreement with FD results is returned from the adjoint MWI with twisted levers (E3).
The relevance of the correct levers is underlined by the iterative convergence behavior that is measured via the discrete global residual to the outer iteration $m \, \widehat{=} \, m^\mathrm{outer}$ based on the L1 norm of the respective local companion, i.e. Eqns. (\ref{equ:adjoint_momentum_balance})-(\ref{equ:adjoint_generic_balance}), viz.
\begin{align}
    R^\mathrm{\varphi, m} = \frac{1}{N} \sum_\mathrm{i=1}^\mathrm{N} |r_\mathrm{i}^\mathrm{\varphi, m} |
    \qquad \qquad \mathrm{with} \qquad \qquad 
    r_\mathrm{i}^\mathrm{\varphi, m}  =  a_\mathrm{ik}^\mathrm{m, \varphi} \varphi_\mathrm{k}^\mathrm{m-1, \varphi} - b_\mathrm{i}^\mathrm{m, \varphi} \, . \label{equ:global_residual}
\end{align}
Here, $a_\mathrm{ik}$ and $b_\mathrm{k}$ denote the system matrix and r.h.s. vector of the equation system belonging to (\ref{equ:semi_discrete_adjoint_momentum}). The residuals are depicted in Fig. \ref{fig:cylinder_adjoint_residuals_drag} over the outer iterations $m$. An adjoint simulation was assumed to converge once the residual of the adjoint velocity -- which controls the shape sensitivity, cf. Eqn. (\ref{equ:sensitivity_derivative}) -- reaches a value of $R^\mathrm{\hat{v}} \leq 10^{-8}$. The residual of the  adjoint velocity (left) -- summed up for both spatial directions --,  adjoint pressure (center) and adjoint concentration equation (right)  are displayed for all MWI configurations E1-E3. In case of E2, the adjoint system converges much slower towards potentially less good sensitivity predictions, cf. Fig. \ref{fig:cylinder_fd_study}. A comparison of E3 and E1 suggests that considering the adjoint momentum sources within an adjoint MWI can lead to a convergence acceleration. For the investigated force functional, the iterative process terminates about $(m^\mathrm{\hat{v},E3} - m^\mathrm{\hat{v},E1}) / m^\mathrm{\hat{v},E1} \cdot 100\% = (20671 - 22281) / 22281 \cdot 100\% \approx 7.2\%$ earlier, cf. Fig. \ref{fig:cylinder_adjoint_residuals_drag} left. 

%
\section{Conclusion}
This note discusses the necessity and implementation of adjoint Momentum-Weighted Interpolation (MWI) strategies. It is outlined that a variety of adjoint momentum sources occur even for simple flows. These sources might feature significant volatile local gradients, which could hamper the convergence. To this end, MWI strategies can be adopted from the primal algorithms and seem viable to improve the robustness of the adjoint solution process. It has been demonstrated that correctly applied face reconstruction methods yield improved sensitivity derivatives and an increased convergence speed. However, the merits of the MWI turn into the opposite and become surprisingly detrimental if the face reconstructions are not carefully implemented.

\end{document}